\documentclass[12pt]{amsart}

\usepackage{amsfonts,amsmath,amsthm}
\usepackage{latexsym}

\newtheorem{theorem}{Theorem}[section]
\newtheorem{corollary}[theorem]{Corollary}

\newtheorem{conjecture}[theorem]{Conjecture}
\newtheorem{proposition}[theorem]{Proposition}

\theoremstyle{definition}

\newcommand{\Z}{\ensuremath{\mathbb{Z}}}
\newcommand{\R}{\ensuremath{\mathbb{R}}}

\newcommand{\Co}{\ensuremath{\mathbb{C}}}
\newcommand{\Spe}{\ensuremath{\mathcal{S}}}

\newcommand{\T}{\ensuremath{\mathcal{T}}}

\def \l {\lambda}
\def \L {\Lambda}

\def \s {\sigma}
\def \S {\Sigma}
\def \C {\Gamma}

\def \W {\Omega}

\def \< {\langle}
\def \> {\rangle}

\begin{document}

\title{Fuglede's conjecture fails in dimension 4}
\author{M\'at\'e~Matolcsi}
\address{ Alfr\'ed R\'enyi Institute of Mathematics,
Hungarian Academy of Sciences POB 127 H-1364 Budapest, Hungary
Tel: (+361) 483-8302, Fax: (+361) 483-8333}
\email{matomate@renyi.hu}
%%% \date{\today}

\begin{abstract}
In this note we give an example of a set $\W\subset \R^4$ such that $L^2(\W )$ admits an orthonormal basis of exponentials $\{\frac{1}{|\W |^{1/2}}e^{2\pi i \< x, \xi \> }\}_{\xi\in\L}$ for some set $\L\subset\R^4$, but which does not tile $\R^4$ by translations. This improves Tao's recent 5-dimensional example, and shows that one direction of Fuglede's  conjecture fails already in dimension 4. Some common properties of 
translational tiles and spectral sets are also proved.
\end{abstract}

\maketitle



\section{Introduction}
Let $\W\subset \R^d$ be a Lebesgue measureable set of finite non-zero measure. A set $\L\subset \R^d$ is said to be a {\it spectrum} of $\W$ if the functions $\{\frac{1}{|\W |^{1/2}}e^{2\pi i \< x, \xi \> }\}_{\xi\in\L}$ form an orthonormal basis of $L^2(\W )$ (here $|\W |$ denotes the Lebesgue measure of $\W$). If such $\L$ exists, $\W$ is said to be {\it spectral}. The set $\W$ is said to be a (translational) {\it tile} if it is possible to tile $\R^d$ with a family of translates $\{t+\W: \ t\in \S\}$ of $\W$ (ignoring overlaps and gaps of measure zero). 

Tiles and spectral sets are connected by Fuglede's famous 

\begin{conjecture}
\cite{fug}A set $\W\subset\R^d$ of finite, non-zero Lebesgue measure is a tile if and only if it is spectral. 
\end{conjecture}
  
Fuglede \cite{fug} proved the conjecture in the special case when the spectrum $\L$ or the translation set $\S$ is assumed to be a lattice.  
The general case of the conjecture was open for nearly 30 years, until last year Tao \cite{tao} showed an example to disprove one direction of the conjecture in 5 and higher dimensions. Namely, he gave an example of a spectral set in $\R^5$ which is not a tile.  The aim of this note is to modify Tao's example and prove that the conjecture fails already in dimension 4. We emphasize that Propositions \ref{prop1} and \ref{prop4} are only new in the generality formulated below, and credit should be given to \cite{tao}, where special cases are proved with the same idea. The essentially new results of the paper are contained in Propositions \ref{prop2}, \ref{prop3} and Theorem \ref{theo1}.  

We will restrict our attention to sets of the form of finite union of unit cubes placed at points with integer coordinates, i.e. 
\begin{equation}\label{kockak}
\W=\cup_{t\in\C}\{t+[0,1)^d\},
\end{equation}
 where $\C\subset\Z^d$ is finite (this means, essentially, that we will be 
 considering
 the tiling and spectral properties of sets in $\Z^d$). Besides giving a counterexample in dimension 4, we also prove some interesting properties  that are shared by spectral sets and tiles of the form \eqref{kockak}.  

To asses our present knowledge concerning this topic it is interesting to note that Fuglede's conjecture was recently proved in 2 dimensions for {\it convex} bodies in \cite{convex}.  
Even for finite union of intervals the conjecture is still open in dimesion 1. Some partial results concerning this case appeared in \cite{lab1}, \cite{lab2}. In arbitirary dimensions certain classes of non-tiling sets are already known to admit no spectrum, either. Results in this direction are contained in \cite{taoref2}, \cite{taoref6} and \cite{taoref10}. For further recent results related to Fuglede's conjecture see \cite{ip98}, \cite{lrw00} \cite{fug01} and \cite{ir02}. The recent survey  \cite{survey} gives an overview not only of the Fuglede conjecture but also of some recent results in the theory of translational tilings.      

Note, also, that one direction of the conjecture is still open in all dimensions, namely we have no examples of tiles which do not admit a spectrum.

\section{Common properties of tiling and spectral sets}

In this section we examine some natural properties of spectral sets and tiles of the form \eqref{kockak}.

We will need the following notations and definitions.

The family of spectral sets (resp. tiles) of the form \eqref{kockak} will
be denoted by $\Spe^d$ (resp. $\T^d$).

A $k\times k$ complex matrix $H$ is called a {\it Hadamard matrix} if all
entries of $H$ have 
absolute value 1, and $HH^\ast =kI$ (where $I$ denotes the identity
matrix). This means that the rows (and also the columns) of $H$ form an
orthogonal basis of $\Co^k$.  A {\it log-Hadamard matrix} is any real square matrix
$(h_{i,j})_{i,j=1}^k$ such that the matrix $(e^{2\pi i h_{i,j}})_{i,j=1}^k$ is
Hadamard. 

For a given finite set $\{t_1, \dots t_k\}\subset \Z^d$ let $T$ denote the
$d\times k$ matrix with $j$th column $t_j$. (With a slight abuse of notation we
will denote
the set itself by $T$ as well.) The set $T=\{t_1, \dots t_k\}$ is
called {\it 
spectral} if there exists a  $k\times d$ real matrix $L$ such that $LT$ is
log-Hadamard (in the usual terminology of harmonic analysis it is also
required that the entries of $L$ belong to the interval $[0,1)$ but the
definition above will be more convenient for present purposes). The (rows of
the) matrix $L$ is the (not necessarily unique) {\it spectrum} of $T$.  The set
$T=\{t_1, \dots t_k\}\subset \Z^d$ is called a {\it tile} if it is possible to
tile $\Z^d$ with a family of translates of $T$.
The family of
spectral sets (resp. tiles) of $\Z^d$ will be denoted by $\Spe_0^d$
(resp. $\T_0^d$).

The set $T=\{t_1, \dots t_k\}\subset \Z^d$ is called {\it m-spectral} for some positive integer $m$ if it
admits a spectrum $L$ with all entries being multiples of $1/m$ (which, essentially,
means that $T$ is spectral in the group $\Z_m^d$). The set $T$ is called an 
{\it m-tile} if it is possible to tile the cube $[0,m)^d$ with translates of $T$ mod $m$. (This means that $T$ tiles the group $\Z_m^d$.)
The family of
$m$-spectral sets (resp. $m$-tiles) of $\Z^d$ will be denoted by $\Spe_m^d$
(resp. $\T_m^d$).

It is obvious that $\Spe_m^d\subset \Spe_0^d$ and $\T_m^d\subset\T_0^d$, and 
it is also clear that if $T\in \T_0^d$ then $T+[0,1)^d\in\T^d$.
It is also well-known that if $T\in \Spe_0^d$ with spectrum $L=\{l_1, \dots
l_k\}$ then $\W=T+[0,1)^d\in\Spe^d$ with spectrum $\Lambda =L+ \Z^d$ (the
orthogonality of the functions $e^{2\pi i\langle \l ,x\rangle}$ is an easy calculation,
while the completeness can be seen by checking the Parseval identity for
functions of the form $e^{2\pi i\< s,x\> }$ and noting that the linear span of
these functions is 
dense in $L^2(\W )$). 

Tao's example \cite{tao} shows that $\Spe_0^d\ne\T_0^d$ and $\Spe^d\ne\T^d$ for $d\ge 5$,  therefore Fuglede's conjecture is false, in general. Nevertheless, it is interesting to see that  $\Spe_0^d$ and $\T_0^d$ share certain properties. Namely, it is reasonable to look for 'natural' properties of $\T_0^d$ and check whether the same holds for $\Spe_0^d$, and vica versa.  
The first result in this direction establishes the 'composition property' for  $\Spe_0^d$ and $\T_0^d$. This property is already used in Tao's example (in less generality, but the idea is the same).

\begin{proposition}\label{prop1}
Assume that the set $T=\{t_1, \dots t_k\}\subset \Z^d$ is an $m$-tile (resp. $m$-spectral)  and the set $S=\{s_1, \dots s_r\}\subset \Z^d$ is an $n$-tile (resp. $n$-spectral). Then the set $\Gamma:=T+mS$ is an $mn$-tile (resp. $mn$-spectral). 
\end{proposition}
\begin{proof}
This property is 'natural' for tiles, therefore we prove only the
spectral part of the proposition.   
  
Let $L:=\{\frac{l_1}{m}, \dots \frac{l_k}{m}\}$ be the spectrum of $T$, and $Q:=\{\frac{q_1}{n}, \dots \frac{q_r}{n}\}$ be the spectrum of $S$. We claim that the spectrum of  $\C$ is $\Lambda=L+\frac{Q}{m}$. 

The cardinality of $\Lambda$ and $\C$ is the same ($mn$), therefore it is enough to check orthogonality of the functions $e^{2\pi i\< \l , x\> }$. 

To see this write $\l_{i_1}=\frac{l_{i_1}}{m}+\frac{q_{i_1}}{nm}$, $\l_{i_2}=\frac{l_{i_2}}{m}+\frac{q_{i_2}}{nm}$, and

\begin{equation}\label{spectcomp}
\sum_{x\in\C}e^{2\pi i \< \l_{i_1}-\l_{i_2},x\> }= \sum_{t\in T}\sum_{s\in S} e^{2\pi i \< \frac{l_{i_1}-l_{i_2}}{m}+\frac{q_{i_1}-q_{i_2}}{nm}, t+ms\> }.  
\end{equation}

Notice that if $q_{i_1}\ne q_{i_2}$ then the second summation gives 0 for all fixed $t\in T$ (because of the spectral property of $Q$ with respect to $S$), therefore the whole sum is 0. 
On the other hand, if $q_{i_1}=q_{i_2}$ then the second summation gives $re^{2\pi i \< \frac{l_{i_1}-l_{i_2}}{m}, t\> }$, and the whole sum equals 0 because of the spectral property of $L$ with respect to $T$.
\end{proof}
         
The next two propositions are not directly used in our 4-dimensional counterexample below, but they are of independent interest.  
First we examine a 'natural' ascendence property of spectral sets.

\begin{proposition}\label{prop2}
Assume that for a given set $T=\{t_1, \dots t_k\}\subset \Z^d$ there exists a $d_1\times d$ matrix $L_1$ with integer entries for which the columns of $L_1T$ form a spectral (resp. $m$-spectral) set $T_1$ in $\Z^{d_1}$. Then $T$ itself is spectral (resp. $m$-spectral).   
\end{proposition}

\begin{proof}
If $L$ denotes a spectrum of $T_1$ then $LL_1$ is clearly a spectrum of $T$, because $(LL_1)T=L(L_1T)$. 

\end{proof}

It may be a little surprising that the same property holds also for tiles. 

\begin{proposition}\label{prop3}
Assume that for a given set $T=\{t_1, \dots t_k\}\subset \Z^d$ there exists a $d_1\times d$ matrix $L_1$ with integer entries for which the columns of $L_1T$ form a tile (resp. $m$-tile) $T_1$ in $\Z^{d_1}$. Then $T$ itself is a tile (resp. $m$-tile).   
\end{proposition}

\begin{proof}
We prove the statement for tiles, and remark that the case of $m$-tiles is settled exactly the same way. 

Note first that the elements of $T_1$ must be different from each other because of the tiling condition.
Let $\S_1$ denote the translation set corresponding to $T_1$, i.e. $\S_1+T_1=\Z^{d_1}$. Define $\S :=L_1^{-1}[\S_1]=\{\s\in\Z^d : \ L_1\s\in \S_1\}$. We claim that $\S$ is a good translation set for $T$. 

To see this, we check first that the translated copies of $T$ are disjoint. Assume therefore that $\s_1+t_i=\s_2+t_j$ for some $\s_1, \s_2\in\S$ and $t_i, t_j\in T$. Applying the transformation $L_1$ we get  
$L_1\s_1+L_1t_i=L_1\s_2+L_1t_j$, and in view of the tiling condition on $T_1$ we conclude that $L_1t_i=L_1t_j$. This implies $t_i=t_j$, and hence $\s_1=\s_2$. 

Next we check that the translates of $T$ cover the whole of $\Z^d$. Let $z\in\Z^d$ arbitrary, and let $z^1=L_1z\in \Z^{d_1}$. Then, by assumption, $z^1= \s^1+t^1$ with some $\s^1\in\S_1$ and $t^1\in T_1$.   
Take the inverse image of $t^1$, i.e. the element $t_i\in T$ for which $L_1t_i=t^1$. Then $\s:=z-t_i$ is contained in $\S$ because $L_1(z-t_i)=z^1-t^1=\s^1\in\S_1$. 
\end{proof}

Let us mention an interesting consequence of this proposition. 

\begin{corollary}
Let $T=\{t_1, \dots t_k\}\subset \Z^d$ be a linearly independent set of vectors (over $\R$). Then $T$ is a tile in $\Z^d$.
\end{corollary}
\begin{proof}
By the assumption of linear independence we can delete $d-k$ rows of the matrix $T$, so that the remaining $k$-dimensional columns $PT=\{Pt_1, \dots Pt_k \}$ span $\R^k$. (The deletion of rows corresponds to multiplication by a $k\times d$ matrix $P$ with entries 1's and 0's.)  In view of Proposition \ref{prop3},  it is enough to show that $PT$ is a tile in $\Z^k$.  

Consider the column vectors $\underline{k}=(0,1,\dots k-1)^*$ (the $*$ in the exponent denotes transposition)  and $\underline{j}:=(\mathrm{det} \ PT) \cdot {PT^*}^{-1}\underline{k}$. Notice that the matrix $(\mathrm{det} \ PT) \cdot {PT^*}^{-1}$ (and  hence the vector $j$) has integer entries, and the 1-dimensional set $\underline{j}^*\cdot PT=(\mathrm{det} PT) \ \underline{k}^*$ obviously tiles $\Z$. Hence, Proposition \ref{prop3} applies, and $PT$ is a tile in $\Z^k$.  
\end{proof}

Finally, we examine a property of 'non-tiling' sets which played an essential
role in the counterexample of \cite{tao} (we give a slightly generalized
version here). 

\begin{proposition}\label{prop4}
Assume that the set $T=\{t_1, \dots t_k\}\subset [0,m)^d$ is {\it not} an $m$-tile. Then for large $n$ the set $S_n:=T+\left ( m\cdot [0,n)^d \right )$ is not a tile. 
\end{proposition}

\begin{proof}
The proof proceeds along the same lines as in \cite{tao}, and we include it
only for completeness.

Assume, for contradiction, that $S_n$ tiles $\Z^d$ with some translation set
$\S$. Take a cube $C_l=[0,l)^d$, where $l$ is much larger than $n$. Let
$\S_l:=\{ \s\in\S \ : (\s + S_n)\cap C_l\ne \emptyset \}$. Note that $\# S_n=kn^d$. We have 
$\#$$\S_l\le\frac{(l+2mn)^d}{kn^d}$, because all $\S_l$-translates of $S_n$ are
contained in the cube $(-mn, l+mn)^d$. 

Let $A$ denote the annulus 
$A:=[-m, mn+m)^d-[m, mn-m)^d$. Then $\# A \  (\approx 4dm(mn)^{d-1})\le
5dm(mn)^{d-1}$, if $n$ is large enough compared to $m$.  Hence, 
$\S_l+A$ cannot cover the cube  $C_{l-m}:=[0,l-m)^d$ because $(\#\S_l)(\# A)\le (l+2mn)^d\left
( \frac{5dm(mn)^{d-1}}{kn^d} \right ) < (l-m)^d$, if the numbers $n,l$ are chosen
so that $n$ is sufficiently large compared to $m$, and $l$ is sufficiently
large compared to $n$. 
 
Take a point $x\in C_{l-m}$ not covered by $\S_l+A$. Consider the cube
$C_m^x:=x+[0,m)^d$. This cube is fully inside $C_l$, therefore if any translate $\s +S_n$ intersects $C_m^x$ then $\s$ necessarliy belongs to $\S_l$. The point $x$ is not covered by the annulus $\s +A$, therefore $C_m^x$ is contained in the cube $\s +[0,mn)^d$. Let $S$ denote the set $T+m\cdot \Z^d$. In view of what has been said, we have $(\s+S_n)\cap C_m^x=(\s+S)\cap C_m^x= (x+[0,m)^d)\cap (\s +T+m\Z^d)$. The mod $m$ reduction of this set is   
exactly the translate $\s+T$ mod $m$. Hence, the tiling of the cube $C_m^x$ 
by 
$\S$-translates of $S_n$  contradicts the assumption that $T$ is not an
$m$-tile. 
\end{proof}

It would be interesting to see whether the corresponding result holds also for spectral sets. 

\section{Counterexample in dimension 4}

We can now show that one direction of Fuglede's conjecture fails already in 
dimension 4. 

\begin{theorem}\label{theo1}
For any $d\ge 4$ there exists a finite union of unit
cubes $\W=\cup_{j=1}^k \{s_j+[0,1)^d\}$  in $\R^d$, such that $\W$ is spectral
but not a tile.
\end{theorem}

\begin{proof}
Take the matrix

$$K:=\left (
\begin{array}{cccccc}
0&0&0&0&0&0\\
0&0&1&1&2&2\\
0&1&0&2&2&1\\
0&1&2&0&1&2\\
0&2&2&1&0&1\\
0&2&1&2&1&0
\end{array}
\right )$$
and observe that $\frac{1}{3}K$ is log-Hadamard (and corresponds to the Hadamard matrix $H$ given in \cite{tao}).

Note that $K$ is of rank 4 mod 3, therefore it is possible to make a mod 3
decomposition $K=L\cdot T$, where $L$ is a $6\times 4$ and $T$ is a $4\times
6$ matrix. We give a specific mod 3 decomposition for convenience:

\begin{eqnarray*}
L:=\left (
\begin{array}{cccc}
0&0&0&0\\
0&1&1&2\\
1&0&2&2\\
1&2&0&1\\
2&2&1&0\\
2&1&2&1
\end{array}
\right )  \ \ \ \ \mathrm{and} \ \ \  
T:=\left (
\begin{array}{cccccc}
0&1&0&0&0&2\\
0&0&1&0&0&2\\
0&0&0&1&0&2\\
0&0&0&0&1&2
\end{array}
\right ).
\end{eqnarray*}

Note, however, that infinitely many decompositions exist, and each
decomposition corresponds to a counterexample as explained below.
(We also remark that the original example of \cite{tao}, although not explicitely
stated in the paper, corresponds to the decomposition $K=K\cdot I$, where $I$
denotes the identity matrix.)

Forgetting about mod 3 calculations and regarding $L$ and $T$ as real matrices we see that $K':=\frac{1}{3}L\cdot T$ is a log-Hadamard
matrix. Therefore the set $T\in\Z^4$ above is 3-spectral, but obviously not a 3-tile because the
number of elements of $T$ (i.e. 6) does not divide $3^4$. 
Therefore, for large $n$, the set $S:=T+3[0,n)^4$ is $3n$-spectral, but not a
tile in view of 
Propositions \ref{prop1} and \ref{prop4}.

This means that the set $\W:=\cup_{s\in S} \{s+[0,1)^4\}$ is spectral  in
$\R^4$ (this fact was already mentioned in the discussion preceding Proposition \ref{prop1}). The fact that $\W$ does not tile $\R^4$ follows as in Proposition
\ref{prop4}. Namely, in a large cube $C_l$ we find a point $x$ which is not
covered by any translate of the annulus $A$, and note that any translate of
$\W$ is either disjoint from the cube $x+[0,3)^4$ or it covers exactly 6 unit
volumes of it. Therefore $x+[0,3)^4$ cannot be covered by translates of $\W$ because of obvious divisibility reasons.

\end{proof}

\end{document}